\RequirePackage{amsmath}
\documentclass{svjour3}                     
\smartqed  
\usepackage{tablefootnote}
\usepackage{multirow}
\usepackage{graphicx}
\usepackage{float}
\usepackage{amsmath}
\usepackage{graphics}
\usepackage{graphicx}
\usepackage{amssymb}

\usepackage{amsthm}
\usepackage{geometry}
\usepackage{float}
\usepackage{color}
\usepackage{picture}
\usepackage[table,xcdraw]{xcolor}
\definecolor{mygray}{gray}{0.8}
\usepackage{multirow}
\usepackage{multicol}
\usepackage{mathtools}

\usepackage[english,ruled,lined,linesnumbered]{algorithm2e}
\usepackage[noend]{algpseudocode}
\SetKwInput{KwData}{input}
\SetKwInput{KwResult}{output}
\makeatletter
\newcommand{\vast}{\bBigg@{4}}
\newcommand{\Vast}{\bBigg@{5}}
\makeatother

\usepackage[authoryear,round]{natbib}
\usepackage{blkarray}
\allowdisplaybreaks
\usepackage{footnote}

\journalname{}
\begin{document}

\title{Heterogeneous Prestressed Precast Beams Multiperiod Production Planning Problem: Modeling and Solution Methods
}


\author{Kennedy Ara\'ujo \and Tib\'erius Bonates \and Bruno Prata \and Anselmo Pitombeira-Neto 
}


\institute{
        K. Ara\'ujo \at
        Institute of Mathematics and Statistics \\
        University of S\~ao Paulo\\
        \email{kennedy94@ime.usp.br}
       \and
        T. Bonates \at 
        Department of Statistics and Applied Mathematics \\
        Federal University of Cear\'a\\
        \email{tb@ufc.br}
        \and
        B. Prata \at
        Department of Industrial Engineering \\
        Federal University of Cear\'a \\
        \email{baprata@ufc.br}
        \and
        A. Pitombeira Neto \at
        Department of Industrial Engineering \\
        Federal University of Cear\'a \\
        \email{anselmo.pitombeira@ufc.br}
}

\date{Received: date / Accepted: date}

\maketitle

\begin{abstract}
 A prestressed precast beam is a type of beam that is stretched with traction elements. A common task in a factory of prestressed precast beams involves fulfilling, within a time horizon, the demand ordered by clients. A typical order includes beams of different lengths and types, with distinct beams potentially requiring different curing periods. We refer to the problem of planning such production as Heterogeneous Prestressed Precast Beams Multiperiod Production Planning (HPPBMPP). We formally define the HPPBMPP, argue its NP-hardness, and introduce four novel integer programming models for its solution and a size reduction heuristic (SRH). We propose six priority rules to produce feasible solutions. We perform computational tests on a set of synthetic instances that are based on data from a real-world scenario and discuss a case study. Our experiments suggest that the models can optimally solve small instances, while the SRH can produce high-quality solutions for most instances.
 \keywords{ prestressed concrete \and precast beams \and modular construction \and integer linear programming \and size reduction heuristics}
\end{abstract}

\section{Introduction}

Unlike conventional precast beams, prestressed precast beams have a different production process, in which they are tensioned using traction elements prior to supporting any actual load. The aim of this process is to improve the resistance and behavior of the beams in service. For their production, a factory uses concrete, along with the traction elements and a set of reusable molds. Traction elements are positioned and tensioned within the molds, after which concrete is cast. This is followed by a curing period, during which the concrete bonds to the traction elements. Those elements are then released and, as their material attempts to resume its original (untensioned) length, the concrete is compressed due to static friction. Prestressed beams are common in factories of civil construction materials, since they are often used in a variety of construction types. They are preferred over steel beams, since concrete has a low price and requires less maintenance when compared to steel. For the purpose of this paper, prestressed precast beams can vary with respect to curing time, length, and the number of traction elements used.

A common task in this type of factory involves fulfilling the demand of a set of clients, within a given time horizon. A typical order includes beams of different lengths and types, with different types of beams potentially requiring different curing times. A mold can be used to produce several beams simultaneously, with the total length of the beams being limited by the mold's capacity. While a given mold can be used to produce different types of beams in different periods, only one type of beam can be produced at a given mold at any given time. The problem of planning such production while minimizing the idle capacity in the molds will be referred to as the Heterogeneous Prestressed Precast Beams Multiperiod Production Planning (HPPBMPP).

To the best of the authors' knowledge, the HPPBMPP is novel, despite its similarity with existing cutting problems. Indeed, we argue that the problem includes a known NP-hard cutting problem as a particular case. The combinatorial nature of the problem makes it hard for managers to generate good schedules in practice, which results in inefficiencies and delays in production. 
The practical importance of the problem also derives from the high-performance, durability, and versatility of prestressed precast beams. Those factors are responsible for the frequent use of such beams in a number of building types and civil structures, ranging from houses and office buildings to bridges and dams. Optimizing the production of prestressed beams has the potential effect of speeding up overall construction time, while improving the usage of molds, allowing factories to accept additional orders due to shorter lead times.


\section{Related Work} \label{relatedwork}

To the best of our knowledge, the HPPBMPP problem has not been previously studied in the revised literature. A special case of the problem has been introduced by ~\cite{de2015integer} and an integer programming model has been proposed for its solution.The authors argued that the problem is closely related to cutting stock and sequencing problems, both of which have been extensively investigated. We bring attention to the following similarities between those problems and the HPPBMPP:

\begin{enumerate}
    \item In the HPPBMPP setting, a mold can  represent a large beam of a certain type that must be cut into smaller pieces, with each piece corresponding to the beams that are produced in the mold. In this interpretation, the leftover part of the large beam corresponds to the mold's unused capacity, rather than actual wasted material;
    \item In HPPBMPP, the production might require several periods before the entire demand has been met, i.e., before all beams have been produced. Producing different beam types may require different curing times. The usage of the molds must be scheduled in such a way as to avoid overlapping (the same mold being used to simultaneously produce different types of beams), while respecting the maximum time allowed, or while minimizing some notion of tardiness.
\end{enumerate}

Cutting stock and sequencing are among the most studied problems in the operations research literature. The one-dimensional cutting problem, in particular, bears close resemblance to the HPPBMPP, in the sense that the production in each mold can be planned (equivalently, the mold can be ``cut'') independently of other molds. In what follows, we highlight some studies that tackle scheduling and cutting problems, and that we consider relevant to our study.

A variety of heuristic methods have been successfully applied to those two classes of problems. \cite{yuen1991heuristics} suggested two heuristics for sequencing cutting patterns in the Australian glass industry and reported substantial savings and low computing times.
\cite{Wascher1996} studied the computational performance of heuristics for the one-dimensional cutting stock problem that work by exploring the neighborhood of an optimal solution to the linear relaxation of a model. The heuristics were reported to find optimal solutions for the majority of the instances tested.
\cite{shahin2004using} presented a genetic algorithm (GA) for solving the one-dimensional cutting stock problem. The authors also studied three real-life scenarios arising from a steel workshop and compared the solutions (cutting schedules) obtained by their algorithm with the actual workshop cutting schedules. \cite{pileggi2005abordagens} presented three heuristic approaches to deal with an integrated pattern generating and sequencing problem. The authors considered the trade-off between the different objective functions involved and compared them in the one-dimensional cutting case.
\cite{benjaoran2005flowshop} proposed a multi-objective flow shop scheduling model for bespoke precast concrete production planning and used a genetic algorithm for its solution. \cite{benjaoran2014three} developed new solution procedures for finding efficient cutting plans while minimizing trim loss and the number of stocks used for the cutting stock problem of construction steel bars. \cite{pitombeira2019} introduced a mixed integer linear programming model to the one-dimensional cutting stock and scheduling problem, and applied successfully a matheuristic based on a fix-and-optimize strategy hybridized with a random local search to solve it.

Studies that are solely based on exact methods as a solution procedure have also been reported.  \cite{arenales2015new} proposed a new mathematical model for the cutting stock/leftover problem (CSLP). Due to the exceedingly large size of the model, the authors proposed to solve its linear relaxation via column generation and to use heuristics for constructing feasible solutions based on the relaxed solution. \cite{braga2016exact} explored an exact and compact assignment formulation for the combined cutting stock and scheduling, along with valid inequalities that are used with a cutting-plane algorithm. 

Another fruitful line of work involves the use of both heuristic and exact methods in a combined solution approach. For instance, \cite{yanasse2007integrated} solved to optimality an integrated problem that involved a cutting stock problem under particular pattern sequencing constraints. Their approach included an integer linear programming (ILP) model, a proposed decomposition scheme to solve the model, a modified subgradient method to solve the dual problem, and several heuristic algorithms. \cite{gramani2006combined} formulated a mixed-integer mathematical model for solving the combined cutting stock and lot-sizing problem in a multi-period planning scenario. The authors proposed a heuristic method based on a shortest path algorithm to minimize trim loss. \cite{nonas2008solving} proposed a new column generating solution procedure for the combined cutting-stock and lot-sizing problem, combined with tree-like and sequential heuristics. \cite{salem2007minimizing} presented three approaches for solving the one-dimensional cutting stock problem: a genetic algorithm, a linear programming model, and an ILP model. The authors studied three real-life case studies from a steel workshop. \cite{arbib2014cutting} proposed an exact ILP formulation for the cutting stock problem with due dates with the aim of minimizing a combination of the number of objects cut and weighted tardiness. The authors developed primal heuristics, upper bounds, and an implicit enumeration scheme. 

The production of precast items has also been previously considered from the optimization viewpoint. As an example, \cite{shih2010optimization} optimized a production project of precast items via a mixed integer linear programming model based on grouping concepts and a recursive procedure. \cite{ko2011precast} approached the problem of scheduling precast production considering six steps: mold assembly, placement of reinforcement and all embedded parts, concrete casting, curing, mold stripping, and product finishing. The authors developed a mathematical model and a multi-objective genetic algorithm to solve it. \cite{augusto2012cost} used genetic algorithms to optimize the design of precast floors. \cite{castilho2012comparative} studied the use of genetic algorithms to minimize the cost of slabs made of
prestressed joists and unialveolar beams. \cite{khalili2013integrated} dealt with the optimization of resources and costs for the precast production of complex configurations by means of a mixed ILP model based on prefabrication configuration and component grouping ideas. \cite{marti2014memetic} analyzed the influence of steel fibers on cost-optimized precast-prestressed concrete road bridges with a double double U-shaped crosssection and isostatic spans using a memetic algorithm with variable depth neighborhood search. \cite{yepes2015cost} proposed a hybrid metaheuristic combining simulated annealing and glowworm swarm optimization to minimize $CO_2$ emissions and cost of the precast bridge production at different stages. \cite{yang2016optimized} made a study in precast production proposing a model for the Flowshop Problem of Multiple Production Lines and developed a genetic algorithm for the problem optimization. The authors identify several objective functions and optimization constraints, although only the optimization objective of makespan minimization was used to simplify the comparisons of the proposed approach. \cite{chen2017optimizing} proposed an ILP model for optimizing precast production planning, allocation of component storage, and transportation, as well as for making timely adjustments for contracted projects, with the aim of minimizing production costs. \cite{wang2018layout} proposed an optimization model for scheduling precast components on pallets during the mold setting process to maximize the average utilization of the pallets and solve it using a constructive heuristic algorithm..

Additionally, some studies have tackled the cutting stock problem considering due dates or multiple periods. For example, \cite{li1996multi} developed heuristics and two two-dimensional cutting stock models with due date and release date constraints, in which meeting orders' due dates are more important than minimizing the waste of materials. \cite{nonaas2000combined} proposed a non-linear optimization model for the combined cutting-stock and lot-sizing problem and suggested several heuristics for finding feasible solutions. \cite{reinertsen2010one} proposed new optimization models for solving the cutting stock problem when orders have due dates. The authors solved the models via column generation, with the corresponding pricing problems solved with shortest path algorithms. \cite{de2015integer} proposed an integer programming model for the multi-period production planning of precast concrete beams. The proposed model, however, handled the simplest case of the HPPBMPP, in which all beams are of the same type, or, equivalently, have unitary curing time.

\section[Problem statement]{Problem statement}
\label{ProblemStatement}

The HPPBMPP consists of planning the usage of the available molds along a given time horizon, i.e., scheduling the beam production in such molds, to cast a demand of prestressed precast concrete beams, possibly of different types, while minimizing the total unused capacity of the molds, i.e. the total idle capacity. In order to formalize the problem, we present the input of the HPPBMPP as follows:

\begin{itemize}
    \item $M$: number of molds in which the beams are produced;
    \item $T$: number of available periods to complete the production;
    \item $C$: number of beam types;
    \item $q_c$: number of distinct lengths of beams of type $c$, with $ c = 1,\ldots,C$;
    \item $l(c,k)$: real numbers corresponding to the actual lengths of beams of type $c$, with $ c = 1,\ldots,C$ and $k = 1, \ldots,q_c$;
    \item $d(c,k)$: demand for beams of type $c$ and length $l(c,k)$, with $c = 1, \ldots,C$ and $k = 1,\ldots,q_c$;
    \item $t_c$: integer number corresponding to the curing time (in terms of periods) of beams of type $c$, for $c=1,\ldots,C$;
    \item $L_m$: real number corresponding to the capacity of the $m$-th mold, with $m = 1,\ldots, M$.
\end{itemize}

\begin{figure}[!ht]
    \centering
    \caption{Example of feasible solution. Three beam types, of various lengths, produced in three molds}
    \includegraphics[scale=0.6]{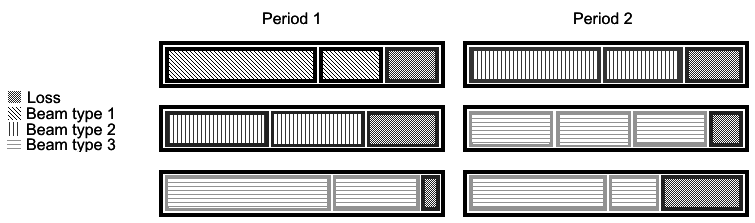}
    \label{exemplo1}
\end{figure}

Each mold can only be used to cast one type of beam at a time. It is possible, however, to simultaneously cast beams of different lengths in the same mold, as long as they are of the same type. The total length of the beams produced during a given period in the $m$-th mold cannot be greater than $L_m$ and the total number of days required to complete the production cannot be greater than $T$. The idle capacity of the $m$-th mold at the $t$-th period, given by $\mathcal{I}(m,t)$, is the difference between the total length of beams produced in the $m$-th mold and its capacity. If the $m$-th mold is not used for beam production in the $t$-th period $\mathcal{I}(m,t)$ is zero. The HPPBMPP output consists of a production plan that minimizes the sum of idle capacities over all molds and periods, i.e., the minimum possible value of $\displaystyle \sum_{t=1}^T \sum_{m=1}^C \mathcal{I}(m,t)$ that can be achieved while fulfilling the demand of beams. An example of a feasible production plan is shown in Figure \ref{exemplo1}.

The HPPBMPP is a combinatorial problem that arises in practical scenarios. Finding an optimal solution to the problem can become a challenging task, as soon as parameters such as the numbers of beam types, lengths and demands increase beyond trivial values. Nevertheless, and despite the similarities between the HPPBMPP and cutting problems, the problem does not precisely fit any existing formulation in combinatorial optimization. Note that, for the purpose of this work, we do not consider neither delivery dates nor stock control in the HPPBMPP and the parameter $T$ is an estimate of the time horizon needed to produce the total demand of beams. We can, however, establish the hardness of the problem:

\begin{proposition}
The HPPBMPP is NP-hard. \label{P1:NPHARD}
\end{proposition}

The assertion in Proposition \ref{P1:NPHARD} can be made due to the fact that HPPBMPP includes, as a particular case, the classical one-dimensional cutting stock problem. Indeed, the case in which there exists only one beam type (i.e., all beams have the same number of cables and the same curing time) turns out to be precisely an instance of the one-dimensional cutting stock problem: the items to be cut correspond to the molds, while the waste of material is equivalent to the unused capacity of each mold. The 1-dimensional cutting stock problem has been known to be NP-hard from the fact that the knapsack problem is reducible to it~\citep{garey1979computers}.

In this section, four models that extend the model proposed by \cite{de2015integer} are described for the case of multiple beam types. In this scenario, different beam types can demand different curing times, unlike the problem treated by \cite{de2015integer}. Moreover, a mold  cannot be used to simultaneously produce beams of different types.

If beams of type $c$ are produced in the $m$-th mold at a given period, it is possible to describe the current state of the mold as a non-negative integer tuple $(a_1,a_2,\ldots,a_{q_c})$, with each $a_k$ ($1 \leq k \leq q_c$) specifying the quantity of beams of length $l(c,k)$ that are currently being produced in the mold. The information of the beam type and the tuple that describes the quantity of each beam length produced --- i.e., the pair $\left(c,(a_1,a_2,\ldots,a_{q_c})\right)$ --- will be called a \emph{pattern}. Naturally, only patterns that do not exceed the $m$-th mold capacity can be produced in that mold. Thus, as a practical matter, we can limit ourselves to taking into consideration only patterns that do not exceed the largest capacity among the molds in the problem's data.

A solution for the HPPBMPP requires fully specifying the pattern that is used in each mold during each of the $T$ periods, with the same pattern being potentially used more than once. The existence of a special pattern $P_0$ is assumed, which is used to denote that a mold is currently being used for the casting of a pattern that began in a previous period and whose production extends at least up to the current period. Since the curing time of each beam type can be different, it is necessary to include constraints in the model that identify the patterns associated with the consecutive periods during which a particular pattern is under production. When our model selects pattern $P_i=\left(c,(a_1,a_2,\ldots,a_{q_c})\right)$, with $c = 1, \ldots, C$, to be initiated in the $m$-th mold at period $t$, it will accordingly select pattern $P_0$ to be used in that mold during the subsequent periods $t+1, \ldots, t+t_c-1$.

For instance, consider that a mold $m$ is used to initiate the production of beams of curing time $3$ at period $5$. Then, the pattern corresponding to the production of those beams must be assigned to period $5$ while $P_0$ must be assigned to the subsequent periods in $m$: $6$ and $7$. This fully describes the state of the mold during periods $5$, $6$, and $7$.

In order to refer to specific information on a given pattern $P_i = \left(\bar{c},(\bar{a}_1,\ldots,\bar{a}_{q_c})\right)$, with $\bar{c} = 1, \ldots, C$, we define the following notation:

\begin{itemize}
    \item $\mathcal{N}_i(c,k)$: number of beams of type $c$ and length $l(c,k)$ that pattern $P_i$ includes. If $c = \bar{c}$, then $\mathcal{N}_i(c,k)=\bar{a}_k$, with $k \in \{1,\ldots,q_c\}$; otherwise, $\mathcal{N}_i(c,k) = 0$, for any $k$.
    \item $u(P_i)$: capacity used by $P_i$, i.e. $\displaystyle u(P_i) = \sum_{k=1}^{q_{\bar{c}}} l(\bar{c},k) \cdot P_i(\bar{c},k)$.
    \item $E_i$: number of periods required to produce the beams in $P_i$. This number equals the quantity of consecutive periods in which $P_i$ remains occupying a mold and is precisely the curing time of beams of type $\bar{c}$, given by $t_{\bar{c}}$.
    \item $F_i^m$: idle capacity of the $m$-th mold when pattern $P_i$ is used in that mold. Note that this quantity depends on the lengths of the beams specified in the pattern, the mold capacity, and the value of $E_i$. $F_i^m$ can be computed as follows: $E_i \cdot (L_m - u(P_i))$. For instance, if the capacity of the $m$-th mold is 10, the capacity used by pattern $P_i$ is $6$, and $E_i = 3$, then we have $F_i^m = 3 \cdot (10 - 6) =  12$.
\end{itemize}

Both $E_i$ e $F_i^m$ can be directly calculated from the problem's data. The value of $F_0^m$, associated to the $P_0$, is defined as zero. However, its is possible to envision variants of the formulation proposed here, in which alternative values for $F_0^m$ are used, depending on the particular objective function to be optimized. 

A remark concerning the use of the $P_0$ pattern is in order. Note that an idle mold (in other words, a mold that is not being used during a specific period) is not assigned the pattern $P_0$. In fact, it has no pattern assigned to it. Moreover, this type of situation is not regarded as a loss. On the other hand, when a mold is used to initiate the production of pattern $P_i$ at period $t$, the subsequent $E_i-1$ periods are assigned $P_0$. This situation results in a total loss of $F_i^m$, which corresponds to the unused capacity of the mold, multiplied by the number of days required for the production of $P_i$.

Given a set of patterns $\{P_1, \ldots, P_r\}$, not including $P_0$, we define the following sets:

\begin{itemize}
    \item $Q(m)$: set containing the indices of the patterns whose capacity does not exceed the capacity of the $m$-th mold: $Q(m) = \{i \in \{1,\ldots,r\} : u(P_i) \leq L_m\}$, for $m=1,\ldots,M$. The same pattern can be used in different molds of potentially distinct lengths. 
    \item $S(j)$: set of indices of the patterns that have curing time $j \in \{1,..., R\}$, with $R$ being the largest curing time of all beam types present in the problem instance.
    \item $Q^*(m)$: set $Q(m)$ including pattern $P_0$, i.e. $Q^*(m) = Q(m) \cup \{0\}$.
\end{itemize}

Our models involve the binary decision variables $x_i^{m,t}$, for $i=1,\ldots,r$, $m=1,\ldots,M$, $t=1,\ldots,T$, each of which is associated with the use of pattern $P_i$ in the $m$-th mold during period $t$, as follows:
\[
x_i^{m,t} = \left\{ \begin{array}{rl}
         1, & \mbox{ if pattern $P_i$ is initiated in the $m$-th mold at period $t$}; \\
         0, & \mbox{ otherwise.}\end{array}\right.
\]

In a scenario of uninterrupted production, exceeding the prescribed demand is usual, although to keep a stock of spare beams can be expensive and limited physically. In addition, in a real-life scenario, it might be desirable to use only patterns that have a minimal percentage of occupation of the molds. If we limit ourselves to using those types of patterns, it may become impossible to satisfy the demands at equality (it could be necessary to use extremely simple patterns to achieve equality). In view of that, the model presented next satisfies the demand with the possibility of surplus. Therefore, the choice of satisfying demands with the possibility of excess seems to be of practical value. 

\subsection{Model for minimizing idle capacity}

We now introduce our main model for the HPPBMPP as follows:

\begin{align}
\textbf{(M1) } \text{minimize} \quad & \nonumber \\
 \quad & \sum_{m=1}^M \; \sum_{i \in Q(m)} \sum_{t=1}^T F_i^m \: x_i^{m,t} \label{FOPPMVPPH}\\
\nonumber \text{subject to} \quad \\
& \sum_{i \in Q^*(m)}\; x_i^{m,t} \leq 1, &\; m=1,\ldots,M, \nonumber \\
&&t=1,\ldots,T  \label{restr1porforma} \\
& \sum_{m=1}^M  \sum_{i \in Q(m)} \sum_{t=1}^{T - E_i + 1} \mathcal{N}_i(c,k) \: x_i^{m,t} \geq d(c,k), & \; k=1,\ldots,q_c, \nonumber\\
&&c=1,\ldots,C  \label{restrdemanda}\\
& (E_i - 1) \; x_i^{m,t} \leq \sum_{\alpha = 1}^{E_i - 1} \; x_0^{m,t + \alpha}, & \; m=1,\ldots,M, \nonumber\\ 
&& \hspace{-4pt} t=1,\ldots,T - E_i + 1, \nonumber\\
&& i \in Q(m) \label{sequenciamento1}\\
&x_0^{m,1} = 0, & m = 1,\ldots, M, \label{sequenciamento3}
\\
&x_0^{m,t} \leq \sum_{\beta = 2}^R \; \sum_{j = \beta}^R \; \sum_{i \: \in \: Q(m) \cap S_j }x_i^{m,t - \beta + 1}, & \; m=1,\ldots,M, \nonumber\\
&& t=2,\ldots,T  \label{sequenciamento2}\\
&x_i^{m,t} \in \{0,1\}, & \; m=1,\ldots,M,  \nonumber \\
&& t=1,\ldots,T, \nonumber\\
&& \label{var1} i \in Q(m) \cup \{0\}. 
\end{align}

The minimization of objective function (\ref{FOPPMVPPH}) has the intent of reducing the total idle capacity of molds and, consequently, concrete waste in used molds. Constraints (\ref{restr1porforma}) ensure that at most one pattern shall be assigned to mold $m$ at period $t$, with the possibility of this pattern being the empty one, $P_0$. Constraint set (\ref{restrdemanda}) requires that all demands be satisfied. Constraints (\ref{sequenciamento1}) force that, if pattern $P_i$ is initiated in mold $m$ at period $t$, then pattern $P_0$ is associated to mold $m$ in the next $E_i - 1$ periods. Note that the right-hand side of the constraint remains unconstrained, in case no patterns is associated to mold $m$ at period $t$. Constraints (\ref{sequenciamento3}) force that pattern $P_0$ is not associated to any mold at the first period. Constraint set (\ref{sequenciamento2}) establish that $P_0$ shall only be used in the $m$-th mold if there is some pattern $P_i$ associated with a previous period in the same mold, such that $P_i$'s production has not yet been completed. Constraints (\ref{var1}) define the domain of the decision variables.

\subsection{Model for minimizing the makespan}

In model (M1), minimizing the objective function (\ref{FOPPMVPPH}) could cause some molds to be unnecessarily filled, particularly in the final periods of production, since the goal was to reduce waste. The next model switches the focus from waste to the time required to fulfill the demand. It makes use of an alternative objective function that captures the number of uninterrupted periods during which at least one mold is used before the production of all the demanded beams is completed. This corresponds to the criterion often used in scheduling problems that measures the time of completion of all jobs, or makespan.

In order to express the minimization of makespan in the model, we introduce another type of decision variable associated with the fact that there is at least one mold used at period $t$:
\[
z_t = \begin{cases}
         1,\mbox{ if at least one mold is used at period $t$, for $t=1,...,T$}; \\
         0,\mbox{ otherwise.}
    \end{cases}.
\]
         
The following model is a specialization of model (M1) and requires the minimization of the makespan:

\begin{align}
\textbf{(M2) } \text{minimize} \quad & \nonumber \\
 \quad & \sum_{t = 1}^T \; z_t \label{makespan}\\
\nonumber \text{subject to} \quad \\
\nonumber &(\ref{restr1porforma})-(\ref{var1})\\
&M \; z_t \geq \sum_{m = 1}^M \left( \sum_{i \in Q^*(m)}\; x_i^{m,t} \right), & t = 1,\ldots, T\label{z_r}\\
\nonumber &\sum_{i \: \in \: Q^*(m) }  x_i^{m,t} \geq \sum_{i \: \in \: Q^*(m) } x_i^{m, t+1}, & m = 1,\ldots, M,  \\
 &&t = 1,\ldots, T-1 \label{continuidade} \\
& z_t \in \{ 0, 1\}, & t = 1,\ldots,T. \label{z_bounds}
\end{align}

Model (M2) includes two additional sets of constraints. Constraint set (\ref{z_r}) ensures that $z_t=1$ when period $t$ is used for the production of any beam, for each period $t$. Note that the periods in which only $P_0$ is used are taken into account by the correspond constraint from (\ref{z_r}), since those patterns correspond to actual production of beams. If no pattern is assigned to any mold during period $t$, then the corresponding variable $z_t$ is not constrained. Since model (M2) minimizes (\ref{makespan}) then $z_t$ will be set to zero whenever (\ref{z_r}) does not impose $z_t=1$.

Minimizing (\ref{makespan}) subject to (\ref{restr1porforma})-(\ref{var1}) and (\ref{z_r})- (\ref{z_bounds}) effectively minimizes the number of days in which production takes place. However, a solution satisfying those constraints might still involve periods of inactivity (i.e., periods in which the production is interrupted), followed by periods of activity. This means that the time of production of the latest beam produced is not necessarily as early as possible. In order to properly capture the makespan of the production plan and accomplish its minimization, we use constraint set (\ref{continuidade}): once the production is interrupted in the $m$-th mold at period $t$, it never resumes. Thus, the complete model minimizes the number of days in which production takes place, while ensuring that those days are contiguous.

A desirable property of model (M2) in the scenario of continuous production is that, once a mold becomes idle, it can be used to start the production of beams to satisfy a demand that was not yet available during the scheduling of the current production plan.

Model (M2) also can be formulated in an alternative way using variable $z$ as an integer variable that defines the makespan in model (aM2):

\begin{align}
\textbf{(aM2) } \text{minimize} \quad & \nonumber \\
 \quad & z \label{makespan2}\\
\nonumber \text{subject to} \quad \\
\nonumber &(\ref{restr1porforma})-(\ref{var1})\\
\nonumber &z \geq \; t \; \sum_{i \in Q^*(m)} \; x_i^{m,t}, & t = 1,\ldots, T,\label{z_uni}\\
& & m = 1,\ldots, M\\
& z \in \{1, \ldots, T \}.
\end{align}

Constraints (\ref{z_uni}) with objective function (\ref{makespan2}) minimization state that $z$ is equal to the index of the last period used to produce beams.

\subsection{Model for minimizing the total completion time}

Although model (M2) generates solutions with less molds unnecessarily filled than model (M1), it still might return solutions with a great amount of unnecessary beams, i.e. more beams than the demand, which have to be stocked. 
Then, one way of avoiding unnecessary beams is minimizing total completion time to achieve the demand. In order to do that, we define model (M3) as follows:

\begin{align}
\textbf{(M3) } \text{minimize} \quad & \nonumber \\
\quad & \sum_{t = 1}^T \sum_{m = 1}^M \Big( \sum_{i \in Q^*(m)}\; x_i^{m,t} \Big) \label{tct}\\
\nonumber \text{subject to} \quad \\
\nonumber &(\ref{restr1porforma})-(\ref{var1}), \; (\ref{continuidade}).
\end{align}

Objective function (\ref{tct}) aims at minimizing the total completion time of beam production.

\subsection{Overview of the models}

Model (M1) involves $\mathcal{O}(M T r)$ decision variables and $\mathcal{O}(q + M T r)$ constraints, with $\displaystyle q=\sum_{c=1}^C q_c$, while model (M2), (aM2) and (M3) have $\mathcal{O}(MTr)$ variables and $\mathcal{O}(q + M T r)$ constraints, as well. We shall not prove this bound formally, as it relies on a straightforward counting argument based on the indices of variables and constraint in models (M1), (M2), (aM2), and (M3). Depending on the total number of possible patterns, there may be an excessive number of variables in all models. In a practical scenario, the unavailability of certain lengths for given beam types might limit the quantity of patterns. Differently, this number can be limited by the exclusive usage of sets of patterns with specific properties. In Subsections \ref{subs:maximal_patterns} and \ref{subs:SRH}, we discuss restrictions on the type of patterns used and why restricting the models in such ways are reasonable approaches.

Models (M1), (M2), and (M3) are linear and have only binary decision variables, a fact that allows for their solution via standard integer linear programming software. A disadvantage of model (aM2) is that it includes one general integer variable. It is interesting to note that all models are also amenable to solution via an iterative scheme of column generation, in which the set of patterns available are generated on demand. This might prove to be useful when dealing with problems that admit very large numbers of patterns.

\subsection{Maximal patterns} \label{subs:maximal_patterns}

A pattern $P_i = (c, (a_1, a_2, \ldots, a_{q_c}))$ is defined as \emph{maximal} with respect to the $m$-th mold if it is not possible to add any beam of type $ c $ to $ P_i $ without violating the capacity of the mold. In our models, we only used variables associated with maximal patterns in their respective molds: that is, a variable $ x_i^{m, t}$ will only exist in the model if $P_i$ is a maximal pattern in mold $m$.

After excluding variables that are associated with non-maximal patterns, there is typically a substantial reduction in the number of variables of all models. This, in turn, can improve their solution times considerably. However, solutions with maximal patterns may produce beam surplus as compared to non-maximal patterns, leading to an increase of the stock size.

We can restrict models (M1), (M2), (aM2) and (M3) to maximal patterns without affecting the optimal values of their objective functions (\ref{FOPPMVPPH}), (\ref{makespan}), (\ref{makespan2}) and (\ref{tct}).

\subsection{Size-reduction heuristic} \label{subs:SRH}

Size-reduction heuristics are solution methods that consist in solving a reduced version of the MILP model in which only a subset of variables is considered. This means that we can drastically reduce the size of the MILP model and depending on the choice of the subset we may be led to promising sub-optimal solutions in shorter execution times and less memory usage. To cite some examples, \cite{fanjul2011size} proposed several size-reduction heuristics for the unrelated parallel machines scheduling problem, reducing the number of machines to only a subset of promising ones taking several criteria into consideration. \cite{fanjul2017models} introduced some matheuristics to the  unrelated parallel machines scheduling problem with additional resources. One of which consists in a size-reduction method named as \textit{job-machine reduction}. Such method involves selecting only variables in which the jobs are associated to the $\ell$ ``best'' machines, otherwise they are removed from the MILP model.

Regarding the HPPBMPP, since the number of maximal patterns still can be too large for state-of-the-art solvers to handle the corresponding MILP model, we can select a subset of maximal patterns to solve the problem, thus not necessarily leading us to an optimal solution, or even a feasible one, for the global problem. Since the number of patterns in the problem is smaller, the number variables and constraints in the MILP model will be smaller.

We define \textit{$q_c$-maximal patterns} as a subset of patterns that are maximal on the shortest mold from an specific instance that covers the largest number of distinct lengths of its beam type. For example, a set $q_c$-maximal patterns in which $q_c = 2$ is a set that has patterns that contain at least 2 beams of distinct lengths. If there is no pattern that covers all beam lengths of a certain type, the set $q_c$-maximal patterns will be composed of by  patterns that covers $q_c-1$ distinct beam lengths, and so on until the set of patterns covers each beam length. Since one characteristic of the problem in practice is that usually there are molds large enough to accommodate a large quantity of beams, it is highly expected that there are patterns that covers all $q_c$ beam lengths for each beam type. 



\section{Priority rules} \label{SecPriority}

In this section we propose six constructive heuristics, which we refer to as \textit{priority rules}, to obtain feasible solutions for the problems under study. Each priority rule that we propose consists in, whenever a mold is freed, selecting a beam type, according to some priority measure regarding the curing time, whose demand has not been attended, and associating it to the current freed mold. Then, we fulfill the current mold with beams of the selected beam type following a second priority measure regarding beam lengths until the demand of the current beam type is achieved or the pattern associated to the current mold is maximal in such mold.

Note that each heuristic described in this section will return solutions that satisfy the demand with no beam surplus, which may lead us to solutions that are composed of patterns that are not necessarily maximal in their respective molds. Regarding this, each of the priority rules that we propose have two phases: the first phase consists in generating a solution producing all demanded beams; the second phase consists in converting each of the patterns used in such solution into maximal patterns. This phase involves filling the patterns with beams of its type from the largest one to the shortest one in matter of length, until each of the generated pattern is maximal in its respective mold.

The priority measures proposed for curing time are:

\begin{itemize}
    \item Shortest curing time first: consists in selecting the beam type with the shortest curing time first among beam types that did not achieved their respective demands;
    \item Longest curing time first: consists in selecting the beam type with the longest curing time first among beam types that did not achieved their respective demands.
\end{itemize}

The priority measures proposed for beam lengths are:

\begin{itemize}
    \item Shortest length first: to select the beam with the shortest length first among beam length from a given type whose demands have not yet been achieved;
    \item Largest length first: to select the beam with the largest length first among beam length from a given type whose demands have not yet been achieved;
    \item Alternate lengths: to select alternately the beam with the shortest length and the beam with the longest length whose demands have not yet been achieved;
\end{itemize}

Based on the measures described above, we name the proposed priority rules as follows:

\begin{itemize}
    \item Shortest curing time shortest length first (SCTSL);
    \item Shortest curing time largest length first (SCTLL);
    \item Shortest curing time alternate length first (SCTAL);
    \item Longest curing time shortest length first (LCTSL);
    \item Longest curing time largest length first (LCTLL);
    \item Longest curing time alternate length first (LCTAL).
\end{itemize}



The complexity of the proposed priority rules is polynomial and it is given by $\displaystyle \mathcal{O}\Big(C \log C + \sum_{c=1}^C q_c \log q_c + M \sum_{c=1}^C \sum_{k = 1}^{q_c} d(c,k) + M T \sum_{c=1}^{C}q_c \Big)$.

\bibliographystyle{spbasic}      
\bibliography{sample}   

\end{document}